\documentclass  [11pt]{article} 
\usepackage                             {latexsym}
\usepackage     [intlimits]             {amsmath} 
\usepackage                             {amssymb}
\usepackage                             {enumerate}
\usepackage                             {theorem}
\theoremstyle                           {break}
\newtheorem     {teo}                   {Theorem}
\newtheorem     {prop}                  {Proposition}
\newtheorem     {lemma}                 {Lemma}

\newtheorem     {definiz}               {Definition}
\newcommand     {\dimo}[1][]            {\noindent\textsc{P r o o f#1.} 
                                        \\ 
                                        \noindent}

\newcommand     {\fine}                 {\newline\noindent%
                                        \textsc{Q.D.E.}\\}
\newcommand     {\ra}                   {\rightarrow}
\newcommand     {\lra}                  {\longrightarrow}
\newcommand     {\demi}                 { \frac {1} {2}  }

\newcommand     {\benum}                {\begin{enumerate}}
\newcommand     {\edenum}               {\end{enumerate}}
\newcommand     {\menuno}               {{-1}}  
\renewcommand   {\phi}                  {\varphi}

\newcommand     {\perdef }              { : \, = }   
\newcommand     {\id  }                 {  \operatorname { Id  }  }

\newcommand     {\restr}[1]             {\phantom{}_{
                                        \text{\raisebox{.4ex}{$|$}}#1}}

\newcommand     {\C}                    { \mathbb { C } }

\newcommand     {\I}                    { \operatorname{\sqrt{-1 }}}

\newcommand     {\di}                   {\, d }

\newcommand     {\est}                  {\raisebox{.34ex}
                                        {$\scriptstyle {\bigwedge}$}}

\newcommand     {\Spin}                 {\operatorname {Spin}}
\newcommand     {\SU}                   {\operatorname {SU}}

\newcommand     {\curv}                 {\operatorname{R}}
\newcommand     {\vol}                  {\operatorname{\mathrm {vol}}}
       
\newcommand     {\meanc}                {\vec{H}}

\renewcommand   {\Re}                   {\operatorname {Re}}

\newcommand     {\dvol}                 {\operatorname { dvol } }
\newcommand     {\grass}                { \operatorname {G} }

\newcommand     {\X}                    { \mathfrak { X} }

\begin{document}

\title{A generalization of Cayley submanifolds}
\author{Alessandro Ghigi}

\maketitle

\begin{abstract}
Given a K\"{a}hler    manifold of complex dimension 4,  we consider  submanifolds of (real) dimension 4, whose K\"{a}hler angles coincide. We call these submanifolds \emph{Cayley}.
We investigate some of their basic properties, and prove that
(a) if the ambient  manifold is  a Calabi-Yau, the minimal Cayley submanifolds are
just the Cayley submanifolds as defined by Harvey and Lawson \cite{hl};
(b) if the ambient is a K\"{a}hler-Einstein manifold of non-zero scalar curvature, then minimal Cayley submanifolds have to be either complex or Lagrangian.
\end {abstract}

\section{Introduction}

Cayley submanifolds were defined by Harvey-Lawson \cite{hl} and by McLean \cite{mclean}, as calibrated submanifolds of $\Spin(7)$-manifolds. Each such manifold $M$ admits a parallel calibration $\Phi \in \est^4 (M)$ whose stabilizer is $\Spin(7)$. It is  called the \emph{Cayley calibration}, due to the link with the octonians, and the corresponding minimal varieties are called \emph{Cayley submanifolds}.

If the ambient manifold is  a Calabi-Yau manifold, the Cayley calibration is not unique. Indeed, given any parallel normalized\footnote {I.e. such that 
\begin{gather}
\frac{\omega^4}{4!} = \vol_g =  \biggl ( \frac{ \I}{  2} \biggr )^4 \Omega \wedge \overline {\Omega} .
\nonumber
\end{gather} 
This amounts to requiring that $||\Omega||^*=1$, where $|| \ ||^*$ denotes the comass norm.} complex volume form $\Omega \in \est^{4,0}(M)$, the form
$$
\Phi_\Omega \perdef \Re \Omega + \frac{\omega^2}{2}
$$
is a calibration whose stabilizer is isomorphic to $\Spin(7)$. Therefore there is  an $\mathrm{S}^1$-family of Cayley calibrations. A submanifold calibrated by any of these forms,  has the property that its K\"{a}hler angles coincide. 
It is therefore natural to consider the submanifolds defined by the latter condition, without any assumption relating to the calibrating forms.  This makes sense whether or not the ambient manifold is Ricci-flat, and gives rise to an interesting family of not necessarily minimal submanifolds, including both the Lagrangian and the complex submanifolds as extreme cases.

In this paper we start collecting some facts from linear algebra, making precise the relation between the K\"{a}hler angles on one side, and the Cayley calibrations on the other. Then we define the submanifolds with equal K\"{a}hler angles, which we call \emph{Cayley}, and prove  a formula \eqref{formola} relating the angle, the mean curvature and the Ricci form of the ambient manifold. 
Finally we apply this formula to the case where the ambient manifold is K\"{a}hler-Einstein, obtaining the following two results:
\begin{teo}\label{I}
Let $(M,J,g)$ be a Calabi-Yau manifold. Then a  Cayley submanifold 
of $M$, is minimal iff it is calibrated by some parallel Cayley calibration.
\end{teo}

\begin{teo}\label{II}
Let $(M,J,g)$ be a K\"{a}hler-Einstein manifold of non-zero scalar curvature.
Then any (connected) minimal Cayley submanifold of $M$ is either complex or (minimal) Lagrangian.
\end {teo}

The first result shows the relation with the theory of Harvey and Lawson.

The last result has been obtained, independently and very recently, also by Isabel Salavessa and Giorgio Valli \cite{sv}, by quite different methods.

\medskip

{\bf Acknowledgments}: 
The author wants to thank Gang Tian for proposing him the subject of this work, and for the constant encouragement. He is also grateful to his advisor, Paolo de Bartolomeis,   and to Claudio Arezzo,  for interesting discussions.

\section {Linear algebra of real 4-planes in $\C^4$}

Let $(V,J,g)$ be a Hermitian vector space of real dimension 8.
Denote by $\omega(X, Y)=g(JX,Y)$  the associated K\"{a}hler form.  
Given a subspace $W \subset V$ we denote by $\pi_{W}$ the orthogonal projection onto $W$, and we put 
$$
B_W \perdef \pi_W \circ J\restr{W}.
$$
$B_W$ is a skew-hermitian operator on $W$ with respect to $g$.
We let $\grass(p, V)$ denote the Grassmannian of \emph{oriented} $p$-planes in $V$. 

Let us recall  an important lemma proved by Harvey and Lawson \cite{kaigai}, applied to our situation.
\begin{lemma}[Canonical form of a 4-plane over $U(4)$] 
Let $(V,J,g) $ be a Hermitian vector space of real dimension $8$. Then,  given $\xi \in \grass(4, V)$, there is a unitary basis $u_1,  u_2, u_3 ,  u_4 $ of $V$ and angles $\theta_1, \theta_2$, with
\begin{equation}
\begin{gathered}
0 \leq \theta_1 \leq \frac{\pi}{2} \\
\theta_{1} \leq \theta_2 \leq \pi 
\end{gathered}
        \label{vinculi}
\end{equation}
such that
\begin{equation}
\begin {split}
  \xi= & u_1 \wedge \bigl ( \cos \theta_1 J u_1 + \sin \theta_1 u_2 \bigr ) \wedge \\
\wedge  &u_{3} \wedge \bigl ( \cos \theta_2 J u_{3} + \sin \theta_2 u_4 \bigr ).
\end{split} \label{formanondigrana}
\end{equation}
Therefore
\begin{gather}
B_\xi = 
\left(         
\begin{array}{cc}
 \begin {array}{cc} 
           0 & -\cos \theta_1 \\
           \cos \theta_1 & 0
       \end{array}               & 0 \\
       0                         &   \begin {array}{cc}
                                        0 & -\cos \theta_2 \\
                                        \cos \theta_2 & 0
                                     \end{array}
\end{array}
\right) \label{repre}
\\
\intertext{and } 
\omega\restr{\xi} = \cos\theta_1 e^{12} + \cos\theta_2 e^{34}.
        \label{HLomega}
\end{gather}
\end{lemma}

The numbers $\theta_1$  and $\theta_2$ are called the \emph{K\"{a}hler angles} of the 4-plane $\xi$. 

\begin{definiz}
$\xi \in \grass (4, V)$ is called a Cayley 4-plane if 
\begin{equation}
  \omega\restr{\xi} = *_{\xi} = \omega\restr{\xi}.
\end{equation}
\end{definiz}
Here $*_\xi$ is the Hodge operator of the metric $g\restr{\xi}$.

\begin{lemma}
An oriented 4-plane $\xi \in \grass(4, V)$ is a Cayley subspace if and only if its K\"{a}hler angles coincide. In this case, putting $\cos\theta_1 = \cos\theta_2 =\lambda \in [0,1]$, we have
\begin{gather}
B_\xi^2 = - \lambda^2 \id \\
( \omega^2 )\restr{\xi}  = 2 \lambda^2 \vol_{g\restr{\xi} } 
\end{gather}
and there is a a positive orthonormal basis $e_1, e_2, e_3, e_4$ of $\xi$ such that
\begin{gather}
\omega\restr{\xi} = \lambda ( e^{12} + e^{34} )  \label{cayomega} \\
B_\xi = 
\left(         
\begin{array}{cc}
 \begin {array}{cc} 
           0 & -\lambda \\
           \lambda & 0
       \end{array}               & 0\\
      0      &   \begin {array}{cc}
                                        0 & -\lambda \\
                                        \lambda & 0
                                     \end{array}
\end{array}
\right) \label{caymatrice}
\end{gather}
\end{lemma}
\dimo
Just apply $*_\xi$ to \eqref{HLomega}.
\fine

A positive orthonormal basis $\{e_1, ..., e_4\}$ in which \eqref{cayomega} hold is called a \emph{Cayley basis}.
If we put 
$$
\X \perdef \{ \xi \in \grass(4,V) : \xi \text{ is Cayley} \},
$$
then we have a well defined function
$$
\lambda : \X \lra [0,1] .
$$
\begin{lemma}
\begin{itemize}
\item [(a)]$\X$ is a closed subset of $\grass(4,V)$ and
$\lambda$ is a continuous function. 
\item [(b)] $\lambda^\menuno (1)$ is the Grassmannian of  complex planes in $(V,J)$, while 
$$
\X_r \perdef \{ \xi \in \X : \lambda (\xi) < 1 \}
$$
consist of totally real subspaces. 
\item [(c)]
$\lambda^\menuno(0) $
is the (oriented) Lagrangian grassmannian, while every $\xi \in \X $ with $ \lambda(\xi) > 0 $ is a symplectic subspace of $(V, \omega)$.
\end{itemize}
\end{lemma}
\dimo
Let us consider the following  subset  of the Stiefel manifold of orthonormal quadruples of vectors in $V$:
\begin{multline}
\mathcal{Y} = \{ (e_1,e_2,e_3,e_4) : 
  \omega(e_1,e_2) = \omega(e_3,e_4) \\
  \omega(e_1,e_3) = \omega(e_1,e_4) 
 =  \omega(e_2,e_3)  = \omega(e_2, e_4)=0 \}
\end{multline}
$\mathcal{Y}$ is a closed subset, and the projection $\pi : \mathcal {Y} \ra \X$ is onto, therefore it is an identification, i.e. $\X$ has the quotient topology. As $\lambda \circ \pi (e_1,e_2,e_3,e_4)=\omega(e_1,e_2)$, $\lambda \circ \pi$ is a continuos function on $\mathcal{Y}$, hence the same is true of $\lambda$. 
The remaining statements are trivial.
\fine

\begin{lemma} 
Let $\xi$ be  a  non-complex, hence  totally real Cayley 4-plane. Given any  Cayley basis $\{e_i\}$ of $\xi$, we put
\begin{equation}
\begin{aligned}
u_1 & = e_1 \phantom { \frac{1} {1 - \lambda^2} }\\
u_3 & = e_3  \phantom { \frac{1} {1 - \lambda^2} }
\end{aligned}
\qquad 
\begin{aligned}
u_2 & = \frac{1} {\sqrt{1 - \lambda^2}} (e_2 - \lambda Je_1 ) \\
u_4 & = \frac{1} {\sqrt{1 - \lambda^2}} (e_4 - \lambda Je_3 ) .
\end{aligned}
\label {baseunitariadiCayley}
\end{equation}
Then $\{u_j\}$ is a unitary basis of $V$ and
\begin{equation}
\begin{split}
\xi  = & u_1 \wedge \Bigl ( \lambda(\xi) Ju_1 + \sqrt{1 - \lambda^2 (\xi) } u_2 \Bigr ) \wedge
\\
\wedge & u_3 \wedge \Bigl( \lambda(\xi) Ju_3 + \sqrt{1 - \lambda^2 (\xi) } u_4 \Bigr).
\label {formanemmenodigrana}
\end{split}
\end{equation}
\end{lemma}
\dimo
A straightforward computation shows that
$$
g (u_i, u_j ) = \delta_{ij} \qquad \omega(u_i, u_j ) =0.
$$
\eqref{formanemmenodigrana} follows immediately from \eqref 
{baseunitariadiCayley}.
\fine

\begin{lemma}   \label {Omegapunto}
\begin{itemize}
\item[(a)]
If $\xi \in \grass(4, V)$ is totally real (i.e. if $\cos \theta_1 \neq 0 \neq \cos \theta_2 $) there exists a unique  normalized
 (4,0)-form $\Omega_\xi$ such that
$$
\Omega_\xi(\xi) > 0.
$$
If we write $\xi$ in the form \eqref {formanondigrana}, then
 $\vec{u} = u_1 \wedge u_2 \wedge u_3 \wedge u_4 $ satisfies $\Omega_\xi (
\vec{u}) = 1$. In particular, two basis $\{u_i\}$ such that \eqref{formanondigrana} hold differ by an element of $\SU(4)$.
\item[(b)]
If $\xi$ is Cayley and totally real, it is calibrated by the Cayley calibration associated to $\Omega_\xi$:
$$
\Phi_\xi = \Re \Omega_\xi + \frac{\omega}{2} \qquad \Phi_\xi (\xi) =1 .
$$
\item[(c)]
If $\xi$ is calibrated by some Cayley calibration, $\Phi_\Omega (\xi)=1$, then it is a Cayley subspace, and $\Omega_\xi = \Omega$.
\end{itemize}
\end{lemma} 
\dimo
From the constraints \eqref{vinculi} descends that
$$
\begin{gathered}
\sin\theta_1 \sin\theta_2 \geq 0 \\
\cos\theta_1 \geq 0.
\end{gathered}
$$
If $\xi$ is totally real, then $\sin\theta_i \neq 0$, and  $\sin\theta_1 \sin\theta_2 >0. $ If we let $\Omega_\xi $ be the unique (4,0)-form such that $\Omega_\xi (\vec{u}) =1$, then $\Omega_\xi (\xi) = \sin\theta_1 \sin\theta_2 >0$. This shows $\Omega_\xi$ only depends on $\xi$ and proves (a).
Using the representation \eqref{formanemmenodigrana} we see that
$$
\Omega_\xi (\xi) = 1 - \lambda^2(\xi) \qquad \omega^2(\xi)= 2 \lambda^2(\xi),
$$
thus proving (b).
On the other hand, using the representation \eqref{formanondigrana} we see that
$$
\Omega_\xi (\xi) = \sin\theta_1 \sin\theta_2
 \qquad \omega^2(\xi)=2 \cos\theta_1 \cos\theta_2.
$$
Therefore, if $\Omega = e^{\I \alpha} \Omega_\xi$, 
\begin{equation*}
\begin{split}
\Phi_\Omega (\xi) 
& = \Re \Bigl ( e^{\I \alpha} \Omega_\xi ( \xi) \Bigr ) + 
        \cos\theta_1 \cos\theta_2 = \\
& = \cos \alpha \sin\theta_1 \sin\theta_2 + \cos\theta_1 \cos=\theta_2
\end{split}
\end{equation*}
and  this can be 1, only if $\cos \alpha=1$ and $\theta_1 = \theta_2$.
\fine

\section{Cayley submanifolds of K\"{a}hler manifolds}

Let $(M, J, g)$ be a K\"{a}hler manifold of complex dimension 4.  We 
consider an oriented submanifold $N \subset M$ of real dimension 4.

We let $*_N$ denote the Hodge operator of the metric $g\restr{N}$.

\begin{definiz}
We call $N$ a Cayley submanifold if the equation
\begin{equation}
\omega\restr{N} = *_N \, \omega\restr{N}
\end{equation}
is satisfied on $N$.
\end {definiz}

This just means that for any point $x$ of $N$, the oriented tangent space $T_x N$ is a Cayley subspace of $T_x M$. 

\medskip

We stress that this definition does NOT agree with the one given by Harvey and Lawson, which makes senses  on any $\Spin(7)$-manifolds and implies that the submanifold is volume-minimizing. The above definition on the contrary makes sense on any K\"{a}hler manifold, and does not imply minimality. Just consider that any Lagrangian submanifold has equal (and zero) K\"{a}hler angles, and is therefore Cayley, according to the above definition.

The relation between this definition and the one of Harvey and Lawson, in the case where the ambient manifold is Calabi-Yau, is the subject of theorem \ref{I}.

\medskip

As the tangent spaces  to $N$ are  Cayley subspaces, 
 if we denote by $B_x$ the endomorphism $\pi \circ J_x \restr {T_x N}$, then
$B_x ^2$ is a multiple of the identity at each point $x$ of $N$. We can define a function $\lambda=\lambda(x) \geq 0$, such that
$$
B_x^2 = - \lambda^2(x) \id.
$$
As $\omega^2 \restr{N}= 2 \lambda^2 \vol$, we deduce that
 $\lambda^2$ is a smooth function on $N$, with values in $[0,1]$.

Given any 4-dimensional submanifold of $M$, not necessarily Cayley, we denote by $N_{r}$  the totally real part of $N$, and by $N_c$ the set of complex  points. If $N$ is Cayley, then 
 $N_r = \{ x \in N : \lambda(x) <1 \}$ and $N_c = \lambda^{-1}(1)$. In particular $N = N_r \sqcup N_c$.  Taking the square root of $\lambda^2$ we deduce that $\lambda$ is a continuos function on $N$,  smooth  on $N_r$, i.e. away from complex points.

On $N_{r}$ is defined a section $\Omega_N$ of $\est^{4,0}(M)\restr{N} = K_M\restr{N}$, determined by the condition that $\Omega_N$ be normalized and satisfy $\Omega_N (T_x N) > 0$ at each point. This is seen applying lemma \ref{Omegapunto}.
\begin{lemma}\label{frames}
Given a Cayley submanifold $N$, near each non-complex point $x $ of $N$,
one can find  a  smooth  Cayley frame $\{e_1, e_2, e_3, e_4\}$, and a smooth unitary frame $\{u_1, u_2, u_3, u_4 \}$ in $TM\restr{N}$ such that 
\begin{equation}
\begin{split}
T_xN  = & u_1 \wedge \Bigl ( \lambda(x) Ju_1 + \sqrt{1 - \lambda^2 (x) } u_2 \Bigr ) \wedge
\\
\wedge & u_3 \wedge \Bigl( \lambda(x) Ju_3 + \sqrt{1 - \lambda^2 (x) } u_4 \Bigr).
\end{split}
\end{equation}
In particular
$$
\Omega_N(u_1, u_2, u_3, u_4) =1,
$$
and therefore
$\Omega_N : N \lra K_M \restr{N} $ is a smooth section.
\end{lemma}
\dimo
For $x \in N_r$, let us consider the endomorphism
$$
j_x = \frac{ B_x}{\lambda(x)}
$$
of $T_xN$. It is a $g\restr{N}$-orthogonal almost complex structure, compatible with the orientation of $N$, smooth on all of $N_r$. Therefore we know that near any $x \in N_r$ we can find a smooth $j$-unitary frame $\{e_1, e_3\}$ in $TN$, i.e. a positive orthonormal frame in $TN$ of the form $\{e_1, j e_1, e_3, je_3 \}$. Putting $e_2 = je_1$, $e_4 = je_3$ we obtain the Cayley frame. Using the formulas \eqref{baseunitariadiCayley} to define $\{u_i\}$ we find a smooth unitary frame of $TM\restr{N}$ with the  desidered properties.
\fine

We let $\nabla$ denote the Levi-Civita connection of the K\"{a}hler metric $g$ on $M$, and $D$ the induced connection on the submanifold $N$. The metric
$g$ being K\"{a}hler, $\nabla$ gives a connection on  $K_M$, and this in turn can be pulled back to a connection on $K_M \restr{N}$. We denote both connections by $\nabla$, too.
Let $\nu_N$ denote the normal bundle to $N \subset M$, $h: TN \otimes TN \ra \nu_N$ the II fundamental form, and $\meanc$ the mean curvature vector.

\begin{prop}
If $N$ is a Cayley submanifold
\begin{gather}
        g(h(X,Y),JZ) - g(h(X,Z),JY) = \bigl ( D_X \omega \bigr ) (Z,Y) 
  \label{eq:simmetrie di h}\\
        \omega(X, \meanc) = \sum_{i=1}^4 g(h(X, e_i) , Je_i)  
\label{simmetriedih2}
\end{gather}
where $X, Y, Z$ are arbitrary vectors tangent to $N$, and $\{e_i\}$ is any orthonormal basis of $TN$.
\end{prop}
\dimo
\begin{equation}
\begin{split}
  g(h(X,Y),JZ) & = g \bigl ( ( \nabla_X Y )^\perp , JZ \bigr ) =
\\
        & = g  (  \nabla_X Y  , JZ ) -  g( D_X Y, JZ) = \\
& = g  (  \nabla_X Y  , JZ ) -  \omega (Z,  D_X Y) ,
\end{split}
\end{equation}
therefore
\begin{multline}
 g(h(X,Y),JZ) - g(h(X,Z)JY) = \\
= g  (  \nabla_X Y  , JZ   )   - g(\nabla_X Z, JY) +
\omega(Y,  D_X Z) - \omega(Z, D_X Y)
\end{multline}
now
\begin{equation}
\begin{split}
g  (  \nabla_X Y  , JZ   )   - g& (\nabla_X Z, JY) = \\
& =X g (Y, JZ) - g (Y, J \nabla_X Z ) - g(\nabla_X Z, JY) = \\
& =X \omega (Z, Y) .
\end{split}
\end{equation}
Therefore 
\begin{equation}
\begin{split}
  g(h(X,Y),JZ) - g      & (h(X,Z)JY) = \\
& = X  \omega (Z, Y) - \omega( D_X Z, Y) - \omega (Z, D_X Y) = \\
& =
\bigl (D_X \omega \bigr ) (Z, Y).
\end{split}
\end{equation}
This proves \eqref{eq:simmetrie di h}. \\
The second formula follows by taking the trace,
$$
\sum_{i=1}^4 \Bigl \{
g( h(e_i, X ) Je_i ) - g (h(e_i, e_i ) , JX) \Bigr \} =
\sum_{i=1}^4 \bigl (D_{e_i} \omega \bigr ) (e_i, X),
$$
and
$$
 \sum_{i=1}^4 \bigl (D_{e_i} \omega \bigr ) (e_i, X) =
-  \di^* \bigl ( \omega\restr{N} \bigr )  (X) .
$$
$N$ being Cayley, the restriction of $\omega$ to $N$ is selfdual (and closed), hence coclosed (with respect to the metric $g\restr{N}$). Therefore $d^* \omega\restr{M} =0$, and 
\begin{equation}
\sum_{i=1}^4 
g( h(e_i, X ) Je_i ) = \sum_{i=1}^4 g (h(e_i, e_i ) , JX) = \omega (X, \meanc).
\end{equation}
\fine

We will now use  a Cayley frame and a unitary basis as in \ref{frames}, defined in some open subset of the totally real part of $N$, to prove some formulas relating $\lambda$ and $\meanc$.

\begin{lemma}
  \begin{equation}
    \Omega_N(u_1, ..., \nabla_X u_k, ..., u_4 ) = \I g( \nabla_X u_k, J u_k ).
  \end{equation}
\end{lemma}
\dimo
$$
\nabla_X u_k = \sum_{j=1}^4 \Bigl \{ g(\nabla_X u_k, u_j ) u_j  +  g(\nabla_X u_k, J u_j ) Ju_j  \Bigr \}
$$
therefore, using the fact that $\Omega_N$ is a complex form, i.e. of type (4,0), we compute
\begin{equation*}
  \begin{split}
\Omega_N(u_1, ..., \nabla_X u_k, & ..., u_4 ) = \\
& = \sum_j g(\nabla_X u_k, u_j )\Omega_N(u_1, ..., u_j, ..., u_4 ) + \\
& +
\sum_j g(\nabla_X u_k, Ju_j )\Omega_N(u_1, ..., Ju_j, ..., u_4 ) =  \\
 =& \,  g(\nabla_X u_k, u_k )\Omega_N(u_1, ..., u_k, ..., u_4 ) + \\
& + g(\nabla_X u_k, Ju_k )\I \Omega_N(u_1, ..., u_k, ..., u_4 ) =  \\
 & = \I  g(\nabla_X u_k, Ju_k ) \Omega_N(u_1, ..., u_k, ..., u_4 ) 
  \end{split}
\end{equation*}
because 
$$
g(\nabla_X u_k, u_k ) = \demi X || u_k ||^2 = 0.
$$
\fine

\begin{lemma}
  \begin{equation}
    ( 1 - \lambda^2 ) \sum_{k=1}^4 g( \nabla_X u_k, J u_k ) = 
\omega (X, \meanc).
  \end{equation}
\end{lemma}
\dimo
Let us use the definition \eqref{baseunitariadiCayley} of  $u_j$:
\begin{align*}
 \begin{split}
\nabla_X & u_2 = \\
&= \biggl (X \frac {1}{\sqrt{1 - \lambda^2}} \biggr ) (e_2 - \lambda Je_1) +
\frac {1}{\sqrt{1 - \lambda^2}} \Bigl ( \nabla_X e_2 - (X\lambda) Je_1 
-\lambda J\nabla_X e_1 \Bigr )
\end{split}
\end{align*}
\begin{align*}
\begin{split}
g( \nabla_X&  u_2,  Ju_2 )  = \\
& = \biggl (X \frac {1}{\sqrt{1 - \lambda^2}}\biggr ) 
\frac {1}{\sqrt{1 - \lambda^2}}
g \Bigl ( e_2 - \lambda Je_1, 
J \bigl (e_2 - \lambda Je_1 \bigr ) \Bigr ) + \\
& + 
\frac{1}{1 - \lambda^2} g\Bigl( \nabla_X e_2 - (X\lambda) Je_1 
-\lambda J\nabla_X e_1, Je_2 + \lambda e_1 \Bigr) = 
\\
& =\, 0 \, + \frac{1}{1 - \lambda^2} g( \nabla_X e_2, Je_2) 
+ \frac{\lambda}{1 - \lambda^2} g(\nabla_X e_2, e_1) +\\
& - \frac{X \lambda}{1 - \lambda^2} 
g( Je_1, Je_2  ) - \frac{X \lambda \,\lambda }{1 - \lambda^2} 
g( Je_1, e_1 ) +\\
&- \frac{\lambda}{1 - \lambda^2} g(\nabla_X e_1, e_2)  +
\frac{\lambda^2}{1-\lambda^2} 
g( \nabla_X e_1, Je_1 )
\end{split}
\end{align*}
therefore
\begin{align*}
 \begin{split}
(1 - \lambda^2) g(\nabla_X u_2, Ju_2 )  
& =
g( \nabla_X e_2, Je_2 ) + \lambda^2 g(\nabla_X e_1, Je_1) + \\
&- \lambda g( \nabla_X e_1, e_2 ) + \lambda g( \nabla_X e_2, e_1 ) = \\
& =
g( \nabla_X e_2, Je_2 ) + \lambda^2 g(\nabla_X e_1, Je_1) + \\
& - g ( \nabla_X e_1, Be_1 ) - g ( \nabla_X e_2, Be_2 )
\end{split} \\
 \begin{split} 
(1 - \lambda^2) \Bigl [ g(\nabla_X u_1 , Ju_1  &) + g( \nabla_X  u_2 , Ju_2 ) \Bigr ]  = \\
& =  g(\nabla_X e_1 ,Je_1 ) + g( \nabla_X e_2, Je_2 ) + \\
& -  g(\nabla_X e_1 ,Be_1 ) - g( \nabla_X e_2, Be_2 ) = \\
& = g(\nabla_X e_1 ,(Je_1)^\perp ) + g( \nabla_X e_2, (Je_2)^\perp ) = \\
& = g( h( X, e_1), Je_1) + g(h(X, e_2) Je_2 ).
\end{split}
\end{align*}
The same computation works for the last two indices, 3 and 4. Summing the two terms and using\eqref{simmetriedih2} one gets
$$
(1 - \lambda^2) \sum_{k=1}^4 g(\nabla_Xu_k, Ju_k) = \sum_i g( h(X, e_i) , Je_i ) = 
\omega(X,\meanc ) .
$$
\fine

\begin{prop}
  \begin{equation}
    \I \bigl ( \nabla_X \Omega_N \bigr ) (\vec{u}) = 
        \frac{i_{\meanc} \omega}{\lambda^2 - 1} (X).
  \end{equation}
\end{prop}
\dimo
By construction  $\Omega_N (\vec{u} ) \equiv 1$. 
\begin{equation*}
0 = X . \Omega_N ( \vec{u}) = \bigl ( \nabla_X \Omega \bigr ) ( \vec{u}) +
\sum_{k=1}^4  \Omega_N ( u_1, ... , \nabla_X u_k, ... , u_4 ).
\end{equation*}
Using the last two lemmas
\begin{equation*}
\I \bigl ( \nabla_X \Omega \bigr ) ( \vec{u}) = 
\sum_{k=1}^4 g(\nabla_Xu_k, Ju_k) =
\frac{ \omega(X,\meanc )} {1 - \lambda^2}.
\end{equation*}
\fine

Let $\rho$ denote the Ricci form of $\omega$ and let us  define $\gamma \in \est^1 (N_{r})$ by 
\begin{equation}
\begin{split}
  \gamma(X) & \perdef \I \bigl ( \nabla_X \Omega_N \bigr )(\vec{u}) = \\ 
  & = \frac{i_{\meanc} \omega}{\lambda^2 -1} (X) =
 \sum_{k=1}^4 g( \nabla_X u_k, J u_k ).
\end{split}
\end{equation}

\begin{teo}\label{III}
  \begin{equation}
    \di \gamma = \rho\restr{N}.
                \label{formola}
  \end{equation}
\end{teo}
\dimo
$$
    X \Bigl ( \bigl ( \nabla_Y \Omega_N \bigr ) (\vec{u}) \Bigr ) =  \Bigl ( \nabla_X \nabla_Y \Omega_N \Bigr ) (\vec{u}) + \Bigl ( \nabla_Y \Omega_N \Bigr ) ( \nabla_X \vec{u} ) .
$$
$(M,g)$ being K\"{a}hler, $\Omega_N \in \est^{4,0}$ implies that $ \nabla_Y \Omega_N \in \est^{4,0}$ too, therefore the same computations as above apply:
\begin{gather*}
\begin{split}
\nabla_Y \Omega_N \bigl (u_1, ...,&  \nabla_X u_k, ..., u_4 \bigr )  = \\
& = 
\sum_{j=1}^4 \Bigl \{ g(\nabla_X u_k, u_j )\Bigl ( \nabla_Y \Omega_N \Bigr )  \bigl (u_1, ..., u_j, ..., u_4 \bigr )  + \\
& +
 g(\nabla_X u_k, Ju_j )\Bigl ( \nabla_Y \Omega_N \Bigr ) \bigl (u_1, ..., Ju_j, ..., u_4 \bigr )  \Bigr \} = \\
 & =  g(\nabla_X u_k, u_k ))\Bigl ( \nabla_Y \Omega_N  \Bigr )\bigl (u_1, ..., u_k, ..., u_4 \bigr )  + \\
& + \I \sum_{j=1}^4 \Bigl \{ g(\nabla_X u_k, Ju_j )\Bigl ( \nabla_Y \Omega_N  \Bigr ) \bigl (u_1, ..., u_j, ..., u_4 \bigr ) \Bigr \} =   \\
 & =  \I g(\nabla_X u_k, Ju_k ) \Bigl ( \nabla_Y \Omega_N  \Bigr ) (\vec{u}),
  \end{split} \\ 
  \Bigl ( \nabla_Y \Omega_N \Bigr ) (\nabla_X \vec{u} ) = \I \Bigl [
  \sum_k g(\nabla_X u_k , Ju_k ) \Bigr ] 
\Bigl ( \nabla_Y \Omega_N \Bigr ) (\vec{u}) = \gamma(X) \gamma(Y).
\end{gather*}
Then applying the usual formula for the differential of a 1-form we find
\begin{gather*}
  \di \gamma (X,Y) = - \I \bigl ( \curv_{XY} \Omega_N \bigr ) (\vec{u}) \\
\intertext{but}
  \curv_{XY} \Omega_N = \I \rho(X,Y) \Omega_N.
\end{gather*}
\fine

Finally we make the following remark.
\begin{prop}
Let $M$ be a K\"{a}hler manifold and let $N, N'$ be two closed Cayley submanifolds in the same homology class,
$$
[N]=[N'] \in H_4 (M, \mathbb{Z}).
$$
Then, if $N$ is Lagrangian, the same is true of $N'$.
\end{prop}
\dimo
It is enough to observe that
$$
|| \lambda_N ||^2 _{L^2 (N)}  = \int_N \lambda^2_N    \dvol = \demi \int_N \omega^2 = \demi < \omega^2 , [N]>
$$
is a topological invariant. If $N$ is Lagrangian, $|| \lambda_{N'} ||^2 _{L^2 }
= || \lambda_{N} ||^2 _{L^2 } =0$, therefore $\lambda_{N'} \equiv 0$ and $N'$ is Lagrangian.
\fine

\section{Minimal Cayley submanifolds in K\"{a}hler-Einstein manifolds}

We now apply the formula \eqref{formola} to the cases where $\rho=s \omega$, i.e. when the ambient manifold is K\"{a}hler-Einstein. This will yield proofs of theorems \ref{I} and \ref{II}.

\dimo[ of theorem \ref{I}]
Let $\Omega \in \Gamma(K_M)$ be a parallel normalized (4,0)-form. On $N$ we can write
$$
\Omega_N = e^{\I \theta} \Omega
$$
for some locally defined real valued function $\theta=\theta(x)$. Then
\begin{align*}
\nabla_X \Omega_N & = \nabla_X \bigl ( e^{\I \theta} \Omega \bigr ) = \\
 & = \I (X \theta ) e^{\I \theta} \Omega = \\
& = \I (X \theta ) \Omega_N. \\
\gamma(X) & = \I \Bigl (\nabla_X \Omega_N \Bigr ) (\vec{u}) = \\
& = - (X \theta)  \Omega_N(\vec{u}) = \\
& = - X \theta
\end{align*}
i.e. $\gamma = - \di \theta$. But $\gamma =0$ because $ \meanc = 0$. Therefore $\theta\equiv \theta_0$ is a constant, and $N$ is calibrated by 
$$
\Phi_0 = \Re \bigl ( e^{\I \theta_0} \Omega \bigr ) + \frac{\omega^2}{2}.
$$
On the other side, if $N$ is calibrated by some parallel Cayley calibration, then it is obviously minimal, and thanks to lemma \ref{Omegapunto} (c), it is Cayley also according to our definition.
\fine

\dimo[ of theorem \ref{II}]
Let $\rho = s \omega$, with $s \neq 0$.\\
If $N \subset M$ is not complex, then $N_{r}$ is not empty. But $\meanc \equiv 0$ implies that $\gamma$, hence $\di \gamma$ vanish identically. Therefore 
$$
\omega\restr{N_r} = \frac{1}{s} \rho\restr{N_r} = \di \gamma = 0
$$
i.e. $N_{r}$ is a Lagrangian submanifold, and $\lambda =0$ on it. This means that  $N_{r} = \lambda^\menuno(0)$ is a closed and open set. Then $N=N_r$, and $N$  is Lagrangian.
\fine


\end{document}